# A New Bi-Objective Model for Resource-Constrained Project Scheduling and Cash Flow Problems with Financial Constraints under Uncertainty: A Case Study

**Seyed-Ali Mirnezami**

*Department of Industrial Engineering, Faculty of Engineering, Shahed University, Tehran, Iran,* mirnezamiali@gmail.com

**Mohammad Ghasemi***

**Corresponding author**, *School of Information Technology, Deakin University, Geelong, 3216, VIC, Australia,* mohammad.ghasemi@research.deakin.edu.au

**Reza Shahabi-Shahmiri**

*School of Industrial Engineering, College of Engineering, University of Tehran, Tehran, Iran,* reza_shahabi011@ut.ac.ir

## Abstract

Owing to the importance of project cash flow, which comprises an entire history of all cash inflows and cash outflows, to economic survival of firms, it is vital to coping with project scheduling issues considering resource constraints in circumstances involving cash flow. Furthermore, since appropriate project management is subject to the innate uncertainties involved in most projects, they are required to be appraised respecting their profound impact. In this paper, a new comprehensive multi-mode multi-objective linear programming model with two conflicting objectives, which are maximizing final cash flow for profit optimization and shortening the duration of project execution, considering improving assumptions, that is, payments delays, project finance constraints, initial capital, different types of interest rates, credit limit to assuage financial distress, credit line usage, is presented in an uncertain environment. Since the model is considered as multi-objective with uncertain parameters, a new extended interval valued fuzzy - Torabi and Hassini (IVF-TH) approach is proposed to tackle the problem. The presented mixed integer linear programming (MILP) model is solved applying CPLEX solver. In addition, a real construction project in oil and gas industry is presented as a case study to illustrate the model applications. Ultimately, for the purpose of assessing the outcomes, a sensitivity analysis is implemented, and the performance of the proposed solution approach is compared to the previous multi-objective optimization methods using both case study and large problem instances.

## 1. Introduction

With respect to the increasing development in a scale of construction projects, there is a desperate need to manage efficiently the complication of modern construction projects (Tavakolan and Nikoukar 2019). The nature of construction projects, in which plenty of difficulties to lessen project duration and expenditures are included, make the project management highly problematic (Cheng et al. 2015). The requirement of declining both the expenses and implementation duration of a project is supposed to be one of the most significant roles of project managers (Zhang and Elmaghraby 2014). Project scheduling problems have been turned out to be the most significant problems, chiefly when project managers are asked to employ resources efficiently with respect to existing resource constraints and with no exceeding project duration (Sahl Abadi et al. 2018). Neglecting to finance while scheduling exerts a deep influence on project cash flow that leads to contractors' encounter with a high rate of project failure (Alavipour and Arditi 2018). Thus, combining financing and scheduling in light of resource constraints has been proven to be crucial for project expenditures as well as project cash flow (Tavakolan and Nikoukar 2019). Unreliable project cash flow anticipation can result in serious failures and place organizations in jeopardy. Therefore, accurate project cash flows are required to be conducted and monitored by project managers over the life of the projects.

Despite the fact that there is an urgent need for project cash flows to be assessed during the project life cycle, project cash flow prognostication in either of project stages has been turned out to be a daunting task (Pate-Cornell et al. 1990). In this regard, several investigations have been conducted from various points of view, including linear and non-linear models, either single or multiple objective functions, different types of resources, and certain or uncertain environments, to cope with the aforementioned problem thus far. Kulejewski et al. (2021) utilized the critical chain method for scheduling the construction of renewable energy facilities. They proposed an approach for sizing the time buffers using deterministic as well as stochastic optimization techniques and considering the optimization criteria that are according to the construction project cash flow analysis. Cevikcan and Kose (2021) presented an optimization method comprising project financing, optimization, and experimental design modules. They employed two mixed-integer linear programming models in a pre-emptive manner to maximize profitability and liquidity. A multi-operator immune genetic algorithm, named MO-IGA, has been presented in Asadujjaman et al. (2022) to cope with resource constrained project scheduling problem (RCPSP) with discounted cash flows. He et al. (2023) proposed a tabu search algorithm for resource-constrained multiproject scheduling problem to minimize the maximal cash flow gap. Zhong and Huang (2023) investigated a discrete-time cash flow problem utilizing dynamic programming principle to solve the problem and considering various interest rates that are controlled using a stationary Markov chain. A centralized resource-constrained multiproject scheduling problem with four different payment patterns was addressed in He et al. (2024). They constructed an optimization model to minimize the maximal cash flow gap and extended a simulated annealing method owing to the NP-hardness of the problem. For the purpose of investigating the impact of integrating supply chain finance (SCF) tools on the cash flow of construction projects, Soliman et al. (2024) compared the traditional financing cycle to an alternative cycle comprising SCF through extending two multi-objective optimization models to maximize the profit and minimize the negative overdraft. In the recent decades, the significance and necessity of taking the innate uncertainties of construction projects into account for cash flow assessment and creation has been fully addressed by a large number of investigations through employing diverse methods such as fuzzy, stochastic, and grey theory (Huang and Zhao, 2014). A chance-constrained programming with a confidence level has been considered by Mirnezami et al. (2023) to tackle multi-skill multi-mode RCPSP under uncertainty. The relationship between cash flow uncertainty and stock price crash risk has been investigated by Wang et al. (2023). Moreover, the most relevant preceding investigations associated with project cash flow are summarily represented in Table 1.

Table 1. Comparisons of the presented model with existent literature

| Researchers | Year | Obj | Financial Requirements | | | | | | | | Resource | | Mode | | | | Uncertainty Approach | |
|---|---|---|---|---|---|---|---|---|---|---|---|---|---|---|---|---|---|---|
| | | | Delayed Payments | Finance Restrictions | Initial Capital | Loans | Daily Expenses | Credit Limit | Interest Rates | Periodic Expenses | Renewable | Non-Renewable | Single | Multi | Scheduling | Cash Flow Modeling | | Case Study |
| San Cristóbal et al. | 2015 | | | | | | | | | ✓ | | | | | | | Gray | ✓ |
| Righetto et al. | 2016 | S | | ✓ | ✓ | ✓ | | | ✓ | ✓ | | | | | | MILP | Robust | ✓ |
| Alavipour & Arditi | 2018 | S | ✓ | ✓ | ✓ | ✓ | | | ✓ | ✓ | | | | | ✓ | LP | | ✓ |
| Alavipour & Arditi | 2019 | M | ✓ | ✓ | | ✓ | | ✓ | ✓ | ✓ | | | | | ✓ | LP | | ✓ |
| Bleyl et al. | 2019 | | | | | | | | ✓ | ✓ | | | | | | | | ✓ |
| Yao & Luo | 2022 | | ✓ | | | | | | ✓ | | | | | | | | Stochastic | |
| Aramesh et al. | 2023 | M | | | | | | | | ✓ | ✓ | ✓ | | ✓ | ✓ | MILP | Interval-valued fuzzy | ✓ |
| Ghasemi et al. | 2023 | M | | | | | | | | ✓ | ✓ | ✓ | | ✓ | ✓ | MILP | Interval uncertainty | ✓ |
| Bruni & Hazir | 2024 | S | ✓ | ✓ | ✓ | ✓ | | | ✓ | | ✓ | | ✓ | | ✓ | MINLP | Robust stochastic | |
| Current study | | M | ✓ | ✓ | ✓ | ✓ | ✓ | ✓ | ✓ | ✓ | ✓ | ✓ | | ✓ | ✓ | MILP | Extended-IVF | ✓ |

As summarized in Table 1, project cash flow problems have predominantly been studied as single-objective problems and only in few investigations multiple objectives have been taken into account. Therefore, there is a significant gap, since despite the fact that the critical objectives of maximizing project profit and minimizing project duration are both essential for contractors and project managers with limited resources, they have been rarely addressed simultaneously. Moreover, while resource limitations have been acknowledged to be as key factors influencing project cash flows, previous studies have not comprehensively incorporated them into their models. Another notable gap lies in managing the inherent uncertainties in project environments. Despite the potential of fuzzy set theory to provide robust solutions, only a limited number of studies have applied this approach to project cash flow problems. Additionally, the execution of activities considering multiple modes has not received sufficient attention in cash flow management problems, even though it aligns better with practical scenarios. Similarly, prior studies have primarily calculated activity expenses on a periodic basis, overlooking the importance of considering both daily and periodic expenses for more accurate financial planning. Finally, while some financial requirements, such as initial capital and credit limits, have a profound impact on cash flow in real-world projects and are required to be addressed concurrently, they have been investigated separately and no study has taken them all into account simultaneously.

This study fills these gaps by introducing a comprehensive multi-objective model that simultaneously considers profit maximization and duration minimization, incorporates resource and financial constraints, applies fuzzy set theory to manage uncertainties, and takes into account both daily and periodic expenses. Furthermore, multi-mode activity execution, which is a more realistic and practical approach for project cash flow management, is addressed in the presented model.

The arrangement of this investigation is illustrated by the following sections. Problem definition is thoroughly addressed in Section 2. Section 3, presents an MILP model with multiple objectives in light of both resources and financial factors. The new proposed approach is introduced in Section 4. Section 5 of this paper is dedicated to demonstrate the proficiency of the presented model utilizing a case study and assess the respective results through performing a sensitivity analysis. Eventually, the conclusion, limitations, and the recommendations for future research of the study are provided in Section 6.

## 2. Problem Description

Project cash flow of the proposed multi-mode multi-objective linear programming model, which aims at maximizing profit while minimizing project duration, is divided into cash inflow and cash outflow. Since activity completion dates rarely align with defined periods, it is often impossible to split an activity's cost across two periods. Therefore, it must be assigned to only one, resulting in an inaccurate project cash flow. For example, in Fig. 1, activity A's expenditures fall entirely within the first period.

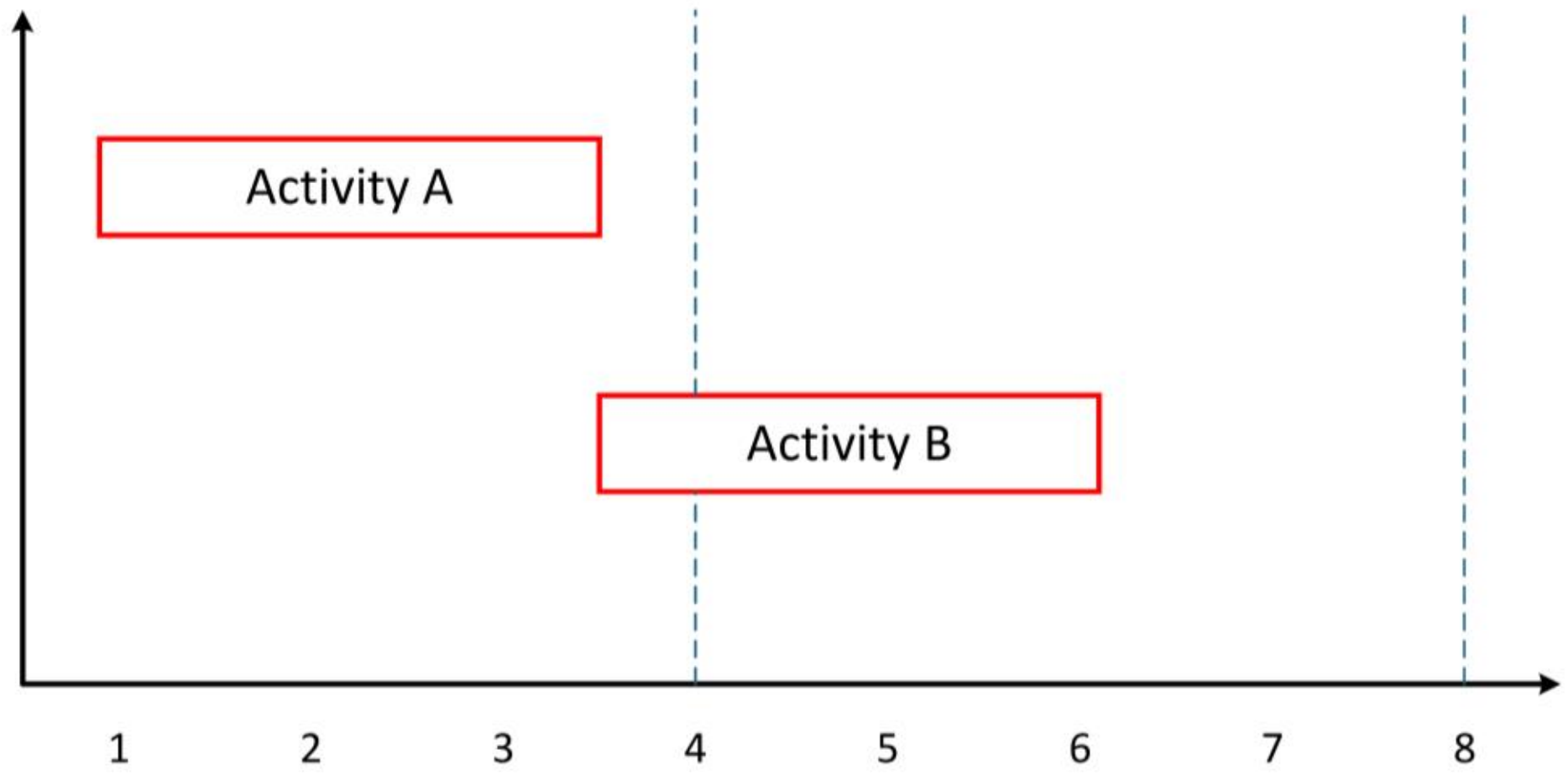

**Fig. 1.** Periodic and daily costs representation

There are circumstances in which start and end dates of activities are not in the same period. For instance, the start date of activity B lies in the first period and its end date is in the second period. Meaning that costs of activity B are divided into the first and second periods. Even though costs of activity B are in two periods, considering these costs periodically leads to dedicating them all to one period, which results in unreliable project cash flow. That is why there is a necessity for activities execution costs, which have been taken into consideration periodically thus far, to be considered both daily and periodically. It is worth noting that contrary to the aforementioned activities costs, which are regarded both daily and periodically, the payments of activities are all devoted to the periods after activities implementation. Therefore, in Fig. 1, the payments of activities A and B are dedicated to the first and second periods respectively.

Moreover, a credit limit is considered to avoid financial distress while project implementation. The required renewable and non-renewable resources for activities accomplishment are restricted to this problem. Considering the aforementioned resource restriction, an optimized project scheduling will be obtained. The initial capital and the long-term loan are supposed to be received only in the first period of project. However, it is assumed that the short-term loans are available to be employed for all periods from the beginning of the project. Furthermore, the basic components of the proposed model are represented in Fig. 2.

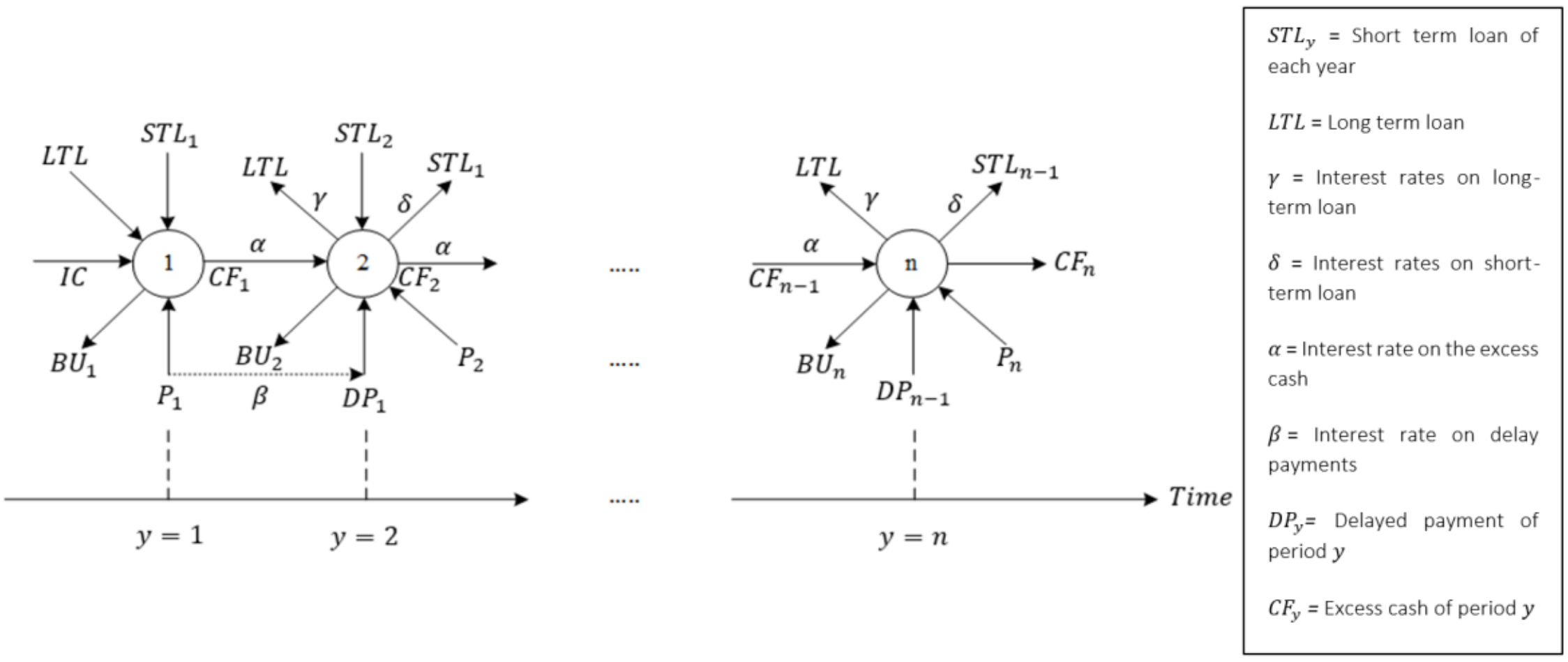

**Fig. 2.** Project cash flow network of the proposed model

Two types of interest rates for these two loans are exerted using:

$$\frac{LTL}{(1+\gamma)^{30}} \tag{1}$$

$$\frac{STL_{y-1}}{(1+\delta)^{30}} \tag{2}$$

Where $LTL$ is the long-term loan, $STL_{y-1}$ indicates the short-term loan of the previous period, and $\gamma$ as well as $\delta$ are interest rates on long-term loan and short-term loan respectively.

Furthermore, two different upper bounds are considered for them, meaning that they are not allowed to be more than a specific amount. The payments of activities execution in various modes are dedicated to periods in which activities are well ordered. Although the proposed model allows for delayed payments if necessary, there is no permission for payments to be delayed more than one period. In the sense that it is possible for payments to be paid incompletely. Thus, the delayed payments of the last period can be received in the next period. It should be noted that the respective interest rate is included in delayed payments and is employed as follows:

$$DP_{y-1} * (1+\beta)^{30} \tag{3}$$

Where $DP_{y-1}$ indicates the delayed payment of previous period and $\beta$ is the interest rate to be employed on delayed payments.

Apart from the introduced interest rates hereinbefore, an interest rate is also considered for excess cash in each period, which is applied by the following:

$$CF_{y-1} * (1+\alpha)^{30} \tag{4}$$

Where $CF_{y-1}$ represents the excess cash of previous period and $\alpha$ is the interest rate on the excess cash.

### *2.1. Assumptions*

The subsequent characteristics are considered in the proposed model as follows. (I) Payments are permitted to be delayed partly or completely. (II) The considered credit lines, that are short and long-term loans, are restricted. (III) Various interest rates for loans, delayed payments, and excess cash are exerted. (IV) The expenditures of activities execution are computed daily and periodically. (V) To assuage financial distress, a credit limit is taken into consideration. (VI) Duration of activities as well as required renewable and non-renewable resources are considered as normalized interval-valued triangular fuzzy (NIVTF) numbers.

### *2.2. Definition of Notations*

Notations applied in this study are determined as follows.

#### *2.2.1. Sets*

| | |
|---|---|
| $i, j = \{1, 2, \ldots, N\}$ | Set of activities |
| $m = \{1, 2, \ldots, M\}$ | Set of activities execution modes |
| $k = \{1, 2, \ldots, K\}$ | Set of renewable resources |
| $l = \{1, 2, \ldots, L\}$ | Set of non-renewable resources |
| $y = \{1, 2, \ldots, Y_n\}$ | Set of long-term periods (monthly) |
| $t, h = \{1, 2, \ldots, T\}$ | Set of short-term periods (daily) |

#### *2.2.2. Parameters*

| | |
|---|---|
| $\tilde{\tilde{d}}_{im}$ | Duration of activity $i$ at mode $m$ |
| $P$ | Precedence relationship between activities |
| $\tilde{\tilde{R}}_{ikm}$ | Renewable resource $k$ usage of activity $i$ at mode $m$ in each day |
| $\tilde{\tilde{W}}_{ilm}$ | Non-renewable resource $l$ usage of activity $i$ at mode $m$ in each day |
| $a_y$ | Lower bound of long-term period $y$ |
| $b_y$ | upper bound of long-term period $y$ |
| $CR_k$ | Cost of using renewable resource $k$ |
| $CW_l$ | Cost of using non-renewable resource $l$ |
| $CC$ | maximum daily cost of using resources |
| $TY_y$ | Interval of period $y$ |
| $PA_{im}$ | Payment for activity $i$ at mode $m$ |
| $\alpha$ | Interest rate on excess cash |
| $\beta$ | Interest rate on delayed payments |
| $\gamma$ | Interest rate on long-term loan |
| $\delta$ | Interest rate on short-term loan |
| $IC$ | Initial capital |
| $\max LTL$ | Maximum long-term loan |
| $\max STL$ | Maximum short-term loan |
| $\min CF$ | Minimum cash flow |

*2.2.3. Decision Variables*

*2.2.3.1. Binary Variables*

| | |
|---|---|
| $X_{imt}$ | Binary variable that equals 1 if activity $i$ is started at mode $m$ on day $t$, otherwise 0 |
| $XP_{imt}$ | Binary variable that equals 1 if activity $i$ is completed at mode $m$ on day $t$, otherwise 0 |
| $XYP_{imyt}$ | Binary variable that equals 1 if activity $i$ is completed at mode $m$ on day $t$ in period $y$, otherwise 0 |

*2.2.3.2. Positive Variables*

| | |
|---|---|
| $BR_{kt}$ | Total renewable resource used on day $t$ |
| $WR_{lt}$ | Total non-renewable resource used on day $t$ |
| $BU_t$ | Cost of using resources on day $t$ |
| $TBU_y$ | Total Cost of using resources in period $y$ |
| $CF_y$ | cash flow in period $y$ |
| $LTL$ | long-term loan |
| $STL_y$ | short-term loan in period $y$ |
| $PA_y$ | Payment of activities in period $y$ |
| $DP_y$ | Delayed Payment in period $y$ |

**3. Mathematical model**

The proposed fuzzy multi-mode multi-objective mixed-integer linear programming model is presented as follows:

$$Min\ Z1 = C_{max} = \sum_{m} \sum_{t} t\ XP_{Nmt} \tag{5}$$

$$Max\ Z2 = CF_{Y_n} \tag{6}$$

$$\sum_m \sum_t X_{imt} = 1 \qquad ;\forall i \tag{7}$$

$$\sum_m \sum_t t\, X_{jmt} \geq \sum_m \sum_t (t + \tilde{\tilde{d}}_{im}) X_{imt} \qquad ;\forall (i,j) \in P \tag{8}$$

$$\sum_i \sum_m \sum_{h=t}^{h=t+\tilde{\tilde{d}}_{im}-1} \tilde{\tilde{R}}_{ikm} X_{imh} \leq BR_{kt} \qquad ;\forall k, \forall t \tag{9}$$

$$\sum_i \sum_m \sum_{h=t}^{h=t+\tilde{\tilde{d}}_{im}-1} \tilde{\tilde{W}}_{ilm} X_{imh} \leq WR_{lt} \qquad ;\forall l, \forall t \tag{10}$$

$$\sum_k CR_k\, BR_{kt} + \sum_l CW_l\, WR_{lt} \leq BU_t \qquad ;\forall t \tag{11}$$

$$BU_t \leq CC \qquad ;\forall t \tag{12}$$

$$\sum_t t\, XP_{imt} = \sum_t (\ t + \tilde{\tilde{d}}_{im}) X_{imt} \qquad ;\forall i, m \tag{13}$$

$$\sum_m \sum_t XP_{imt} = 1 \qquad ;\forall i \tag{14}$$

$$\sum_y XYP_{imyt} = XP_{imt} \qquad ;\forall i, m, t \tag{15}$$

$$TY_{y-1} XYP_{imyt} \leq t\, XP_{imt} \qquad ;\forall i, m, y, t \tag{16}$$

$$t\, XYP_{imyt} \leq TY_y\, XYP_{imyt} \qquad ;\forall i, m, y, t \tag{17}$$

$$\sum_i \sum_m \sum_t PA_{im}\ XYP_{imyt} - PA_y \leq DP_y \qquad ;\forall y \tag{18}$$

$$\sum_{t \geq a_y}^{b_y} BU_t = TBU_y \qquad ;\forall y \qquad (19)$$

$$CF_y = IC + STL_y + LTL + PA_y - TBU_y \qquad ;\forall y = 1 \qquad (20)$$

$$CF_y = STL_y + CF_{y-1} * (1+\alpha)^{30} + PA_y + DP_{y-1} * (1+\beta)^{30} - TBU_y - \frac{LTL}{(1+\gamma)^{30}} - \frac{STL_{y-1}}{(1+\delta)^{30}} \qquad ;\forall y \geq 2 \qquad (21)$$

$$LTL \leq \max LTL \qquad (22)$$

$$STL_y \leq \max STL \qquad ;\forall y \qquad (23)$$

$$CF_y \geq \min CF \qquad ;\forall y \qquad (24)$$

$$X_{imt}, XY_{imyt}, XP_{imt}, XYP_{imyt} \in \{0.1\}, BR_{kt}, WR_{lt}, BU_t, TBU_y, CF_y, LTL, STL_y, PA_y, DP_y \in Z^+ \qquad ;\forall i, m, y, t \qquad (25)$$

The objective function (5) minimizes project duration through minimizing the completion time of last activity while the second objective, that is objective function (6), maximizes the project profit through maximizing the cash flow of the last period. Constraint (7) ensures the implementation of each activity in one time and mode. Constraint (8) represents the precedence relationships between activities. According to this constraint, activity $j$, which is a successor activity, is only permitted to be commenced after accomplishment of activity $i$, which is a predecessor activity. Constraints (9) and (10) determine the required renewable and non-renewable resources in each day respectively. Constraint (11) indicates the total cost of using both kind of resources in each day and the next constraint, constraint (12), is an upper bound for this cost. Constraint (14) ensures the completion date of each activity, which is determined in constraint (13), in one time and mode. Constraints (15) – (17) establish relationship between the periods and the completion date of activities through devoting the periods to them. The amount of delayed payments in each period is determined in constraint (18). The total daily costs of utilizing resources are summarized for each period in constraint (19). Apart from the denoted cash flow of the first period by constraint (20), from the second period on, constraint (21) is in charge of cash flow representation. Two upper bounds are taken into account for the required credit lines, namely long and short-term loans, in constraints (22) and (23) respectively. Furthermore, to avoid financial distress, constraint (24) considers a lower bound for cash flow periodically. Eventually, both the binary variables and the positive variables of the model are defined by constraint (25).

## 4. Proposed Solution Approach

Owing to the innate uncertainties in construction projects, which lead to unhealthy project cash flow, and their profound impacts upon cash flow prognostication, taking them into account while estimating the parameters has been proven to be essential. Activity duration and resource requirement are two critical uncertain parameters that can directly impact on financial planning and project feasibility, particularly in cash flow problem. With respect to the superiority of fuzzy approach over stochastic approach respecting the abovementioned issues, fuzzy mathematical programming is applied in this paper to tackle uncertainty (Moreno and Blanco 2018).

In comparison with generalized fuzzy numbers, interval-valued fuzzy (IVF) numbers as a particular form of them, which are first presented by Gorzalczany (1987), have been turned out to be more flexible (Chatterjee and Kar 2016) and intelligent (Chen and Chen 2008). Furthermore, in numerous cases, the boundaries of generalized fuzzy numbers effectively address ambiguity. As a result, Interval-Valued Fuzzy (IVF) numbers are utilized to provide a more robust representation. Particularly in cash flow optimization, where uncertain resource requirements and fluctuating activity durations directly influence financial considerations, IVF numbers present a more robust representation of these variations. The ability of IVF numbers to model imprecise boundaries makes them appropriate for capturing the variability in cash inflows and outflows, ensuring a more resilient decision-making framework. In this study, normalized interval-valued triangular fuzzy (NIVTF) numbers are taken into consideration. Based on the study carried out by Yao and Lin (2002), interval-valued triangular fuzzy numbers in a minimization problem can be expressed as follows: $\tilde{\tilde{A}} = [\tilde{A}^L, \tilde{A}^U] = [[a_o^L, a_m^L, a_p^L; \tilde{h}_{\tilde{A}}^L], [a_o^U, a_m^U, a_p^U; \tilde{h}_{\tilde{A}}^U]]$, where $a_o^U < a_o^L < a_m^L < a_m^U < a_p^L < a_p^U$, meaning that $\tilde{A}^L \subset \tilde{A}^U$ and membership functions are denoted as $\mu_{\tilde{A}}^{L(U)} = \tilde{h}_{\tilde{A}}^{L(U)}$. IVF can be turned into NIVTF if and only if $a_m^L = a_m^U$ and $\tilde{h}_{\tilde{A}}^L = \tilde{h}_{\tilde{A}}^U = 1$. Fig. 3 represents a NIVTF $\tilde{\tilde{A}}$ (Stanujkic, 2015).

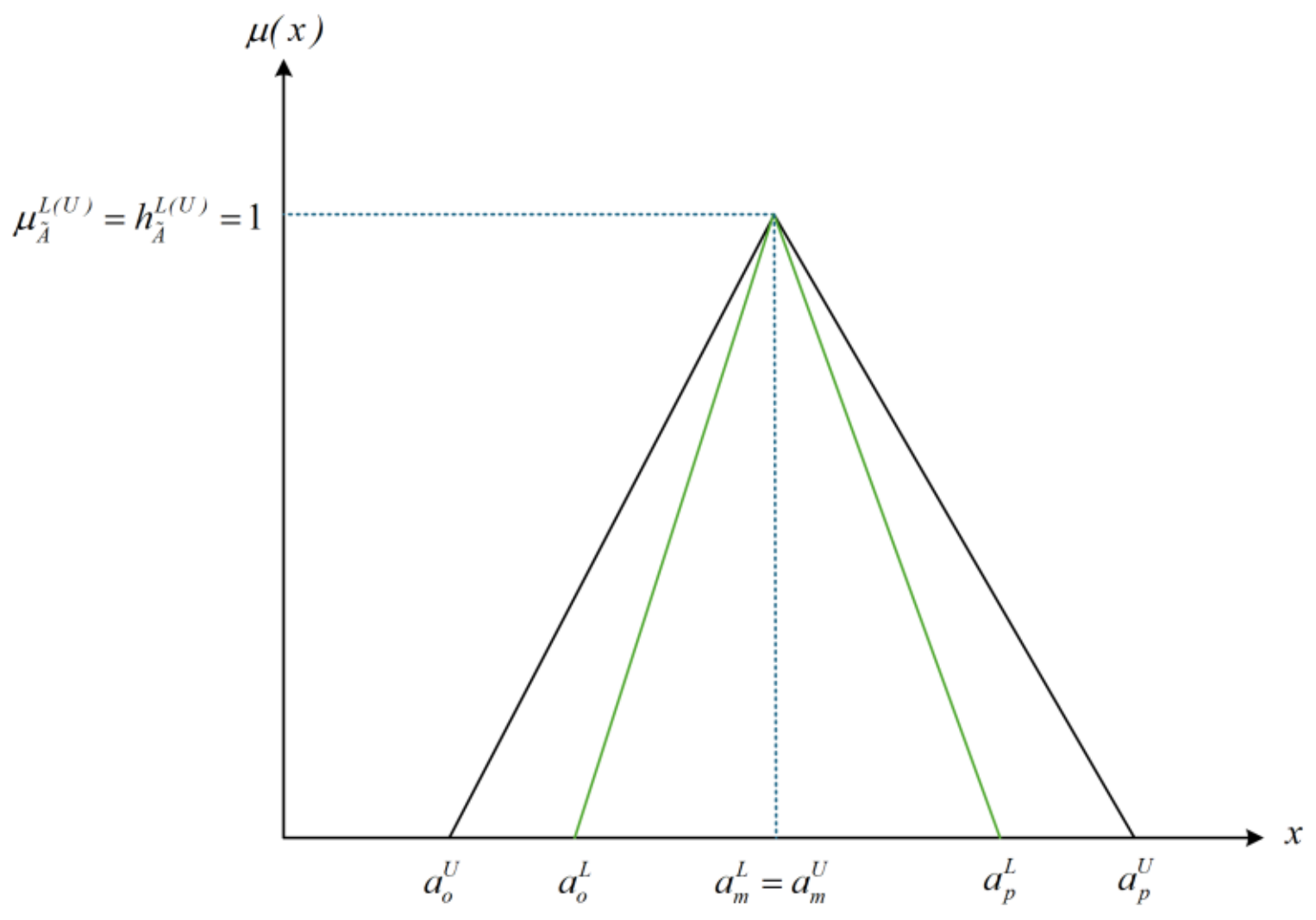

**Fig. 3.** A normalized interval-valued triangular fuzzy (NIVTF) number $\tilde{\tilde{A}}$ (Stanujkic, 2015)

### *4.1. A New Extended IVF-TH Approach*

Since the proposed mixed-integer linear programming model is regarded in a fuzzy environment, an efficient approach is exerted in accordance with Jiménez et al. (2007) and Parra et al. (2005) to convert the model into the equivalent auxiliary crisp model. According to Jiménez et al. (2007), both the expected value and the expected interval of the NIVTF numbers will be computed. A fuzzy mathematical programming model, whose all parameters are as NIVTF numbers, is considered. The following model can be regarded as the equivalent crisp $\alpha$-parametric model of the abovementioned model (Jiménez et al. 2007).

$$\min\ [EV(\tilde{\tilde{c}}^T)^{L(U)}]x$$

Subject to:

$$((1\text{-}\alpha)E_2^{a_i^{L(U)}} + \alpha E_1^{a_i^{L(U)}})x \geq \alpha E_2^{b_i^{L(U)}} + (1-\alpha)E_1^{b_i^{L(U)}} \qquad i = 1,2,\dots,l$$

$$((1\text{-}\tfrac{\alpha}{2})E_2^{a_i^{L(U)}} + \tfrac{\alpha}{2}E_1^{a_i^{L(U)}})x \geq \tfrac{\alpha}{2}E_2^{b_i^{L(U)}} + \left(1-\tfrac{\alpha}{2}\right)E_1^{b_i^{L(U)}} \qquad i = l+1,\dots,m$$

$$(\tfrac{\alpha}{2}E_2^{a_i^{L(U)}} + (1-\tfrac{\alpha}{2})E_1^{a_i^{L(U)}})x \leq \left(1-\tfrac{\alpha}{2}\right)E_2^{b_i^{L(U)}} + \tfrac{\alpha}{2}E_1^{b_i^{L(U)}} \qquad i = l+1,\dots,m$$

$$x \geq 0 \qquad (26)$$

Where $[EV(\tilde{\tilde{c}}^T)^{L(U)}] = \frac{(c^o)^{L(U)}+2(c^m)^{L(U)}+(c^p)^{L(U)}}{4}$, $E_1^{a_i^{L(U)}} = \frac{(a_i^o)^{L(U)}+(a_i^m)^{L(U)}}{2}$, $E_1^{b_i^{L(U)}} = \frac{(b_i^o)^{L(U)}+(b_i^m)^{L(U)}}{2}$, $E_2^{a_i^{L(U)}} = \frac{(a_i^m)^{L(U)}+(a_i^p)^{L(U)}}{2}$, and $E_2^{b_i^{L(U)}} = \frac{(b_i^m)^{L(U)}+(b_i^p)^{L(U)}}{2}$.

Since the proposed model is considered as a multi-objective model, there is a need for its optimal solutions to be acquired through the respective single-objective form (Hadi-Vencheh and Mohamadghasemi 2015). Several solution approaches such as the methods of Torabi and Hassini (2008), and Alavidoost et al. (2016), have been employed for the purpose of obtaining appropriate solutions for multi-objective optimization models. Among previous investigations, Torabi and Hassini (2008) 's method is adopted to be applied. As such, $\alpha$-positive ($OF_i^{\alpha-PIS}$) and α-negative ($OF_i^{\alpha-NIS}$) ideal solutions initially are determined. Afterward, the membership degree of each objective function can be computed using the following membership functions.

$$\mu_i(x) = \begin{cases} 1 & OF_i > OF_i^{\alpha-PIS} \\ \dfrac{f(x) - OF_i^{\alpha-NIS}}{OF_i^{\alpha-PIS} - OF_i^{\alpha-NIS}} & OF_i^{\alpha-NIS} \leq OF \leq OF_i^{\alpha-PIS} \\ 0 & OF_i < OF_i^{\alpha-NIS} \end{cases} \tag{27}$$

$$\mu_i(x) = \begin{cases} 1 & OF_i < OF_i^{\alpha-PIS} \\ \dfrac{OF_i^{\alpha-NIS} - f(x)}{OF_i^{\alpha-NIS} - OF_i^{\alpha-PIS}} & OF_i^{\alpha-PIS} \leq OF \leq OF_i^{\alpha-NIS} \\ 0 & OF_i > OF_i^{\alpha-NIS} \end{cases} \tag{28}$$

Ultimately, with respect to acquired membership functions from the previous step, the following single-objective model can be employed. Thus, appropriate optimal solutions of the multi-objective crisp model will be found.

Where both $\gamma$ and $\theta_i$ indicate the weights, which are determined by decision maker, such that $\gamma \in [0, 1]$, $\theta_i \in [0, 1]$ and $\sum_{i=1}^{2} \theta_i = 1$. Fig. 4 shows the overall procedure of the solution approach.

In this paper, a new extended IVF-TH approach is introduced according to Torabi and Hassini (2008) to deal with the multi-mode multi-objective linear programming model for cash flow problem in a fuzzy environment in which Torabi and Hassini's (2008) method and the normalized triangular interval-valued fuzzy numbers are considered simultaneously. The following model can be regarded as the equivalent crisp $\alpha$-parametric model of the abovementioned fuzzy model (Jiménez et al. 2007) with the same objective functions. For project managers, accurately defining a deterministic value for two critical project parameters, namely activity durations and required resources, can be particularly challenging in real-world projects due to inherent uncertainties and complexities. It may not be possible for deterministic values to be determined. Therefore, project managers can provide subjective judgment of the potential values for these parameters. However, given the inherent uncertainty in their judgments, they often lack complete confidence in these estimations. As a result, they may present these values as an interval range, considering the possible variations and providing a more flexible approach for decision-making. Due to the lack of sufficient historical data for the mentioned parameters, they cannot be considered as stochastic parameters.

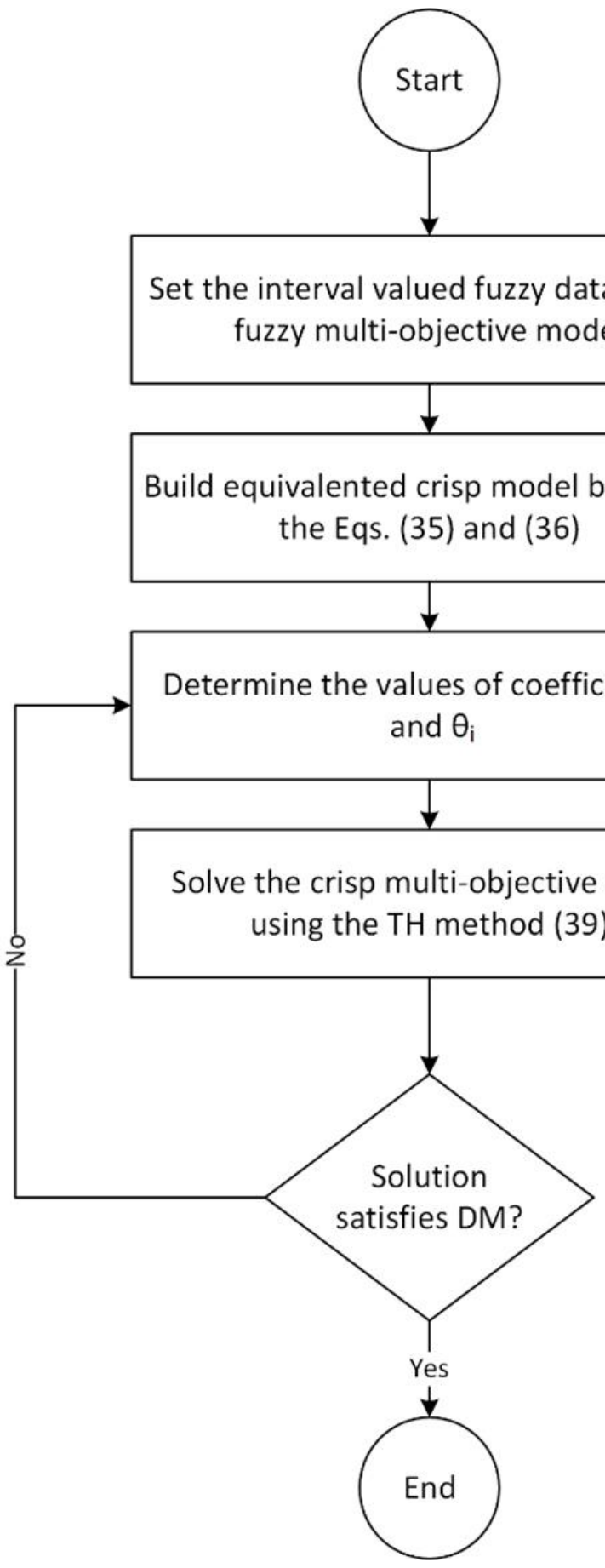

**Fig. 4.** An overview of proposed solution approach

With respect to the inherent ambiguity in the boundaries of generalized fuzzy numbers, activity durations ($\tilde{\tilde{d}}_{im}$), required renewable resources ($\tilde{\tilde{R}}_{ikm}$), and required non-renewable resources ($\tilde{\tilde{W}}_{ilm}$) are considered as NIVTF numbers in this study. NIVTF numbers enable decision-makers to deal with lack of information and uncertainties by relying on intuition and personal experiences rather than solely concrete facts or evidence. Consequently, since these parameters only appear in constraints (8)–(10) and (13), the respective constraints are reformulated as follows:

$$\sum_{m}\sum_{t} t\, X_{jmt} \geq \sum_{m}\sum_{t}(t + \alpha E_2^{d_{im}^L} + (1-\alpha)E_1^{d_{im}^L}))X_{imt} \qquad ;\forall (i,j) \in P \qquad (30)$$

$$\sum_{i}\sum_{m}\sum_{h=t}^{h=t+(\frac{d_{im}^p+2d_{im}^m+d_{im}^o}{4})-1} ((1-\alpha)E_2^{R_{ikm}^L} + \alpha E_1^{R_{ikm}^L})\, X_{imh} \leq BR_k \qquad ;\forall k, \forall t \qquad (31)$$

$$\sum_{i}\sum_{m}\sum_{h=t}^{h=t+(\frac{d_{im}^p+2d_{im}^m+d_{im}^o}{4})-1} ((1-\alpha)E_2^{W_{ilm}^L} + \alpha E_1^{W_{ilm}^L})\, X_{imh} \leq WR \qquad ;\forall l, \forall t \qquad (32)$$

$$\sum_t t\, XP_{imt} \geq \sum_t \left( t + \frac{\alpha}{2} E_2^{d_{im}^L} + (1 - \frac{\alpha}{2}) E_1^{d_{im}^L} \right) X_{imt} \qquad ; \forall i, m \qquad (33)$$

$$\sum_t t\, XP_{imt} \leq \sum_t \left( t + (1 - \frac{\alpha}{2}) E_2^{d_{im}^L} + \frac{\alpha}{2} E_1^{d_{im}^L} \right) X_{imt} \qquad ; \forall i, m \qquad (34)$$

Constraints (6), (11), (12), and (14) – (25)

Taking the following constraints into account and applying the proposed extended IVF-TH approach, the model can be solved and both $Z_i^L; i\epsilon\{1,2\}$ and $Z_i^U; i\epsilon\{1,2\}$ will be computed.

$$Z_1^U \leq Z_1^L \qquad (35)$$

$$Z_2^U \geq Z_2^L \qquad (36)$$

## 5. Performance Appraisement of the Presented Model

### *5.1. Case Study*

A real-world gas refinery construction project is presented as a case study in this paper to illustrate the proposed mixed-integer linear programming model applications. With respect to the importance of construction projects, especially in oil and gas industry, experienced corporations are required to be examined. Among these corporations in this area, MAPNA company is adopted to be investigated. Due to the significance of combinations of rich gas such as ethane, propane, and butane, the Kangan Petro Refining company was established. One of phases produced 81 million cube meters rich gas daily. The 12$^{th}$ phase was located in the east southern part of South Pars Gas Field in the border line with Qatar. In this phase, rich gas was used without recovering ethane. On 14$^{th}$ February 2016 a contract between Kangan Petro Refining company as owner and MAPNA company entered into force. Regardless of the number of both completed and ongoing projects, it should be noted that one of them is represented in Fig. 5 and is investigated in this paper.

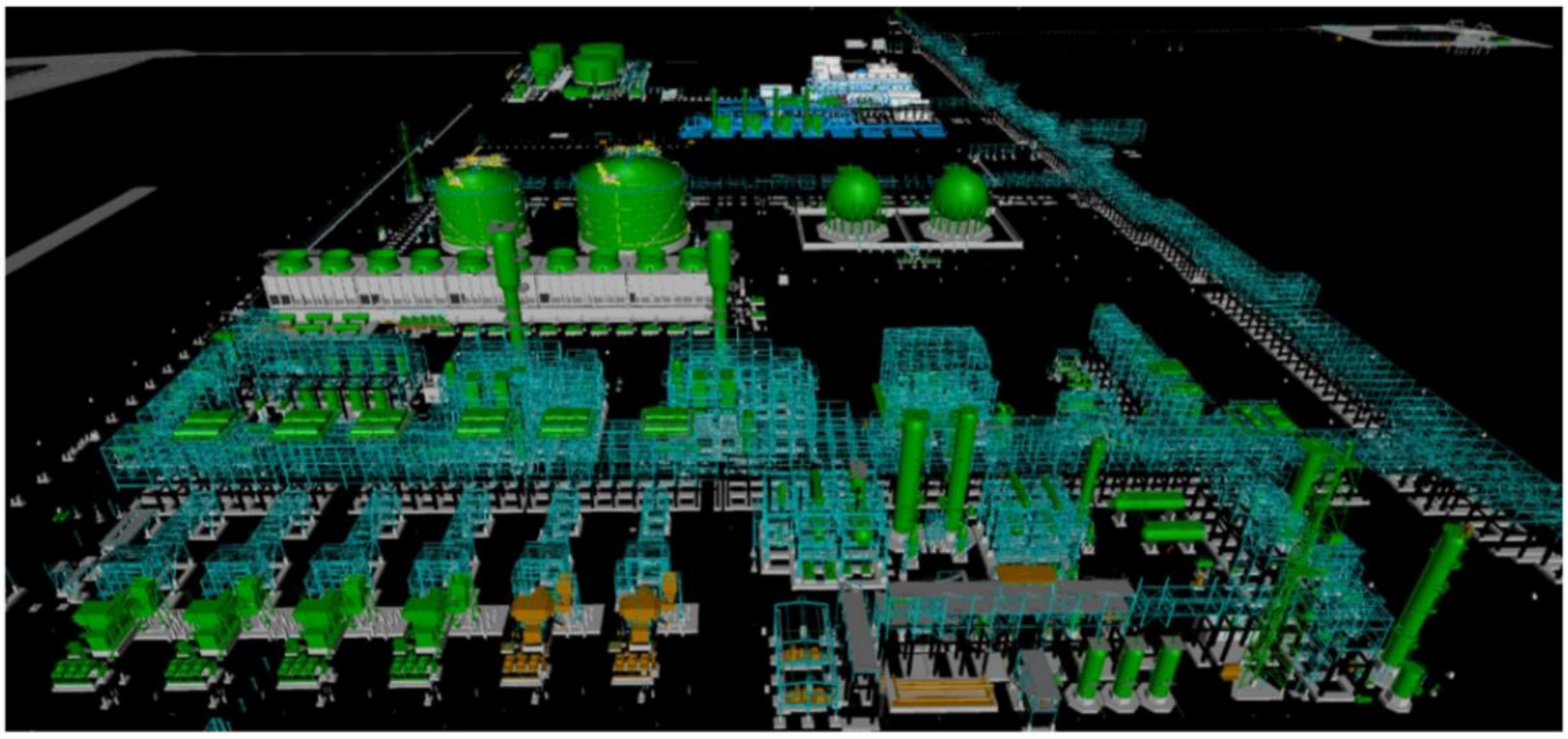

**Fig. 5.** An illustration of the project

### *5.2. Implementation*

As it was mentioned before, although several projects have been implemented by the company, one of them with twenty-two activities is taken into account in this paper. This project, that is steam generation, consists of different activities, including construction, precommissioning, start-up, and performance test, for boilers B, C, D, and E. However, only piping of boiler E is considered in the current study. For better understanding, the work breakdown structure (WBS) of the project is depicted in Fig. 6.

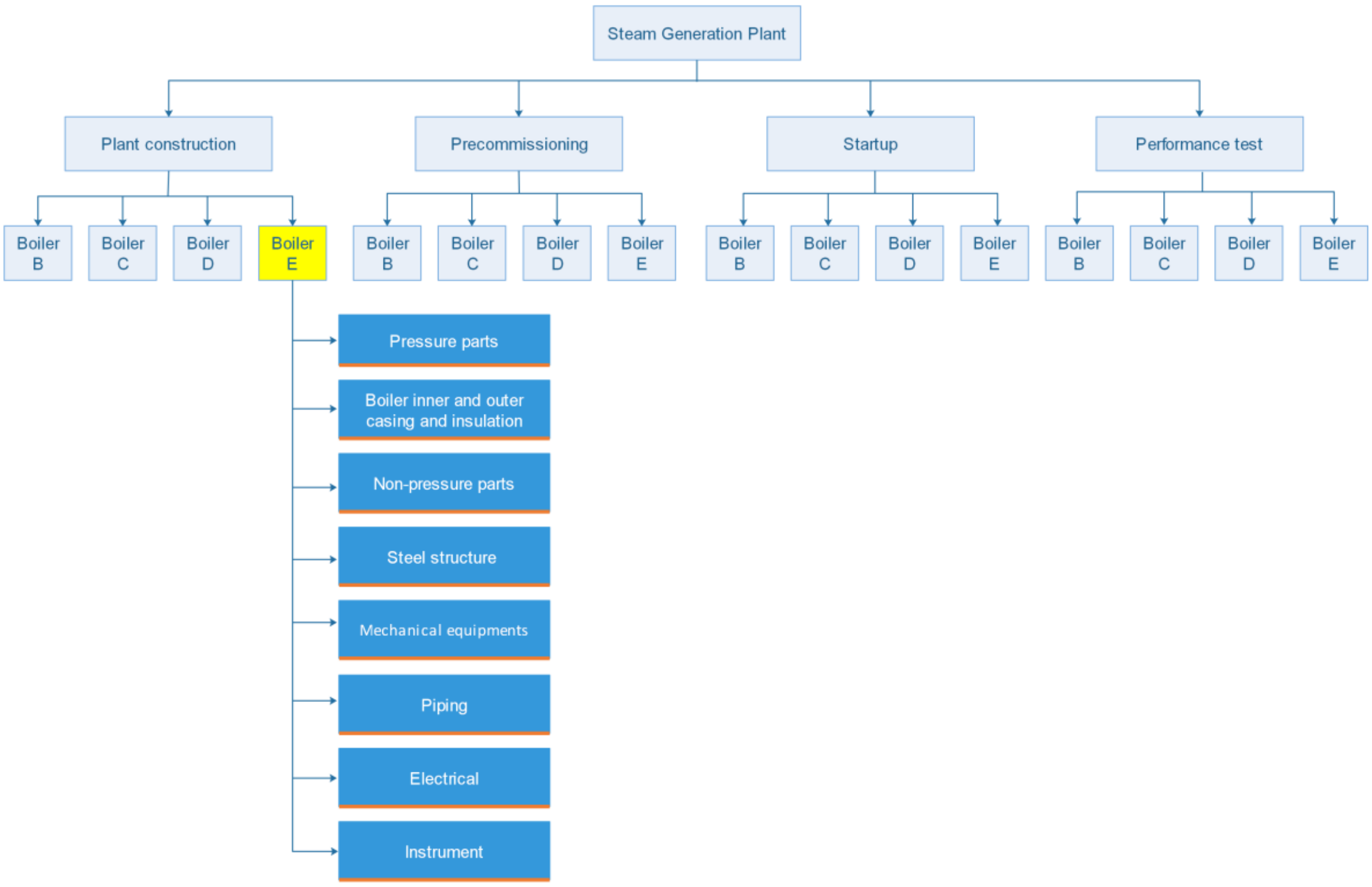

**Fig. 6.** WBS of project

The utilized data in this paper are all in accordance with the adopted project of the corporation documentation. These data, including activity durations, and payments are listed in Table 2.

Table 2. Comparisons of the presented model with existent literature

| Activities | Predecessors | Duration (days) | | Payment ($) | |
|---|---|---|---|---|---|
| | | Mode 1 | Mode 2 | Mode 1 | Mode 2 |
| Start | - | - | - | - | - |
| Material handling of piping - Boiler E | Start | 1 | 1 | 2925 | 4355 |
| Positioning and leveling carbon steel piping - Boiler E | Material handling of piping - Boiler E | 20 | 18 | 4615 | 7085 |
| Fitting-up and welding carbon steel piping - Boiler E | Positioning and leveling carbon steel piping - Boiler E | 10 | 7 | 7800 | 9815 |
| NDT and PWHT of carbon steel of piping - Boiler E | Fitting-up and welding carbon steel piping - Boiler E | 5 | 4 | 7735 | 9620 |
| Positioning stainless steel piping - Boiler E | Positioning and leveling carbon steel piping - Boiler E | 10 | 9 | 9165 | 10465 |
| Leveling stainless steel piping - Boiler E | Positioning stainless steel piping - Boiler E | 5 | 3 | 9295 | 10660 |
| Fitting-up and welding stainless steel piping - Boiler E | leveling stainless steel piping - Boiler E | 10 | 8 | 11115 | 12220 |
| NDT of stainless steel piping - Boiler E | Fitting-up and welding stainless steel piping - Boiler E | 5 | 3 | 7020 | 8450 |
| Cutting galvanized piping - Boiler E | Fitting-up and welding stainless steel piping - Boiler E | 15 | 12 | 6955 | 8060 |
| Threading galvanized piping - Boiler E | Cutting galvanized piping - Boiler E | 10 | 6 | 8905 | 10010 |
| Installation of galvanized piping - Boiler E | Threading galvanized piping - Boiler E | 10 | 8 | 8775 | 10465 |
| Installation of accessories for galvanized piping - Boiler E | Installation of galvanized piping - Boiler E | 5 | 3 | 9035 | 10855 |
| Installation of valve and accessories for piping - Boiler E | Installation of galvanized piping - Boiler E | 15 | 12 | 9815 | 11570 |
| Fabrication of pipe support - Boiler E | Positioning and leveling carbon steel piping - Boiler E | 5 | 2 | 11245 | 12935 |
| Installation of pipe support - Boiler E | Fabrication of pipe support - Boiler E | 5 | 3 | 11635 | 13520 |
| Welding of pipe support - Boiler E | Installation of pipe support - Boiler E | 20 | 17 | 10595 | 12545 |
| NDT of pipe support - Boiler E | Welding of pipe support - Boiler E | 10 | 8 | 8320 | 10205 |
| Touch-up of pipe support - Boiler E | Welding of pipe support - Boiler E | 10 | 9 | 10270 | 12675 |
| Sandblasting piping - Boiler E | Installation of galvanized piping - Boiler E | 5 | 4 | 8645 | 10985 |
| Applying primer to piping - Boiler E | Sandblasting piping - Boiler E | 10 | 7 | 8255 | 10530 |
| Intermediate coating - Boiler E | Sandblasting piping - Boiler E | 5 | 4 | 7085 | 9620 |
| Final coat and touch-up of piping - Boiler E | NDT of stainless steel piping - Boiler E, Installation of accessories for galvanized piping - Boiler E, Installation of valve and accessories for piping - Boiler E, NDT of pipe support - Boiler E, Sandblasting piping - Boiler E, Applying primer to piping - Boiler E, and Intermediate coating-Boiler E | 1 | 1 | 4550 | 6240 |
| End | Final coat and touch-up of piping - Boiler E | - | - | - | - |

The first and second columns in Table 2 are allocated for the representation of project activity names and the relevant predecessors respectively. Depending on how the activities are supposed to be executed, various modes can be determined. Therefore, duration of activities as well as the respective payments in different modes are represented in the following columns.

Applying case study data as input of model and implementing the proposed mixed-integer linear programming model through the CPLEX Solver in GAMS 24.1.2, which is run on a PC with 2.30 GHz CPU and 6 Gb RAM, with the corresponding CPU time of 964 s, all the optimal solutions are found. The respective optimal solutions are listed in Table 3 and Table 4. The optimal solution of objective functions, which are minimizing the project makespan and maximizing the profit, are represented in Table 3. As such, the least project duration and the most profit respecting the considered constraints in the problem are 114 days and $4,914,100, respectively. A point to note is that, the presented case study is applied in both uncertain and crisp environments, and the project durations are obtained 67 and 114 in these circumstances respectively, which indicates the importance of taking uncertainties into account.

Table 3. The derived optimal solutions

| Objective functions | Optimal solutions ($Z_1^*, Z_2^*$) |
|---|---|
| OF 1 (Project duration) | 114.000 days |
| OF 2 (Profit/ Final cash flow) | $ 4.9141E+6 |

Table 4. The results stemmed from the proposed model

| Periods variables | Period 1 | Period 2 | Period 3 | Period 4 |
|---|---|---|---|---|
| Cash flow ($) | 890,000.000 | 1474,232.399 | 2809,530.777 | 4,914,108.476 |
| Total cost of using resources ($) | 210,000.000 | 253,350.000 | 401,100.000 | 327,100.000 |
| Short-term loan ($) | 100,000.000 | 100,000.000 | 0 | 100,000.000 |
| Payments ($) | 0 | 0 | 0 | 79,495.000 |
| Delayed payments ($) | 19,890.000 | 62,010.000 | 56,355.000 | 0 |

Project duration is divided into four periods and the optimal solutions for other determined decision variables of the problem are denoted in Table 4. In this regard, the cash flow of the project in which all cash inflows and outflows are included, is calculated in each period. As it is mentioned before, not only the daily expenditures but also the periodic expenditures of activities execution are taken into account. It is worth noting, although both types of expenditures are considered, one of which is periodic expenditure is represented in this table. According to the calculations, as there is no need for long-term loan, meaning that its value is zero, short-term loans are only denoted herein. Since one of objective functions is maximizing the profit, interest rate is considered for delayed payments, and as noted in the assumption, there is no limitation for delayed payments, the payments of activities are postponed to the last period. That is why there are no payments in the first three periods and they are all delayed in these periods. Given the considered limitations, duration of activities, and available resources, Gantt chart of the project, in which *x*- and *y*-axis are project duration and activities respectively, is illustrated in Fig. 7. The start time, finish time, execution mode, and duration of activities are presented in this figure. In accordance with this Gantt chart, project duration is 114 days. A point to note is that, as the first and the last activities are dummy activities and their duration is zero, the respective start and finish time are the same.

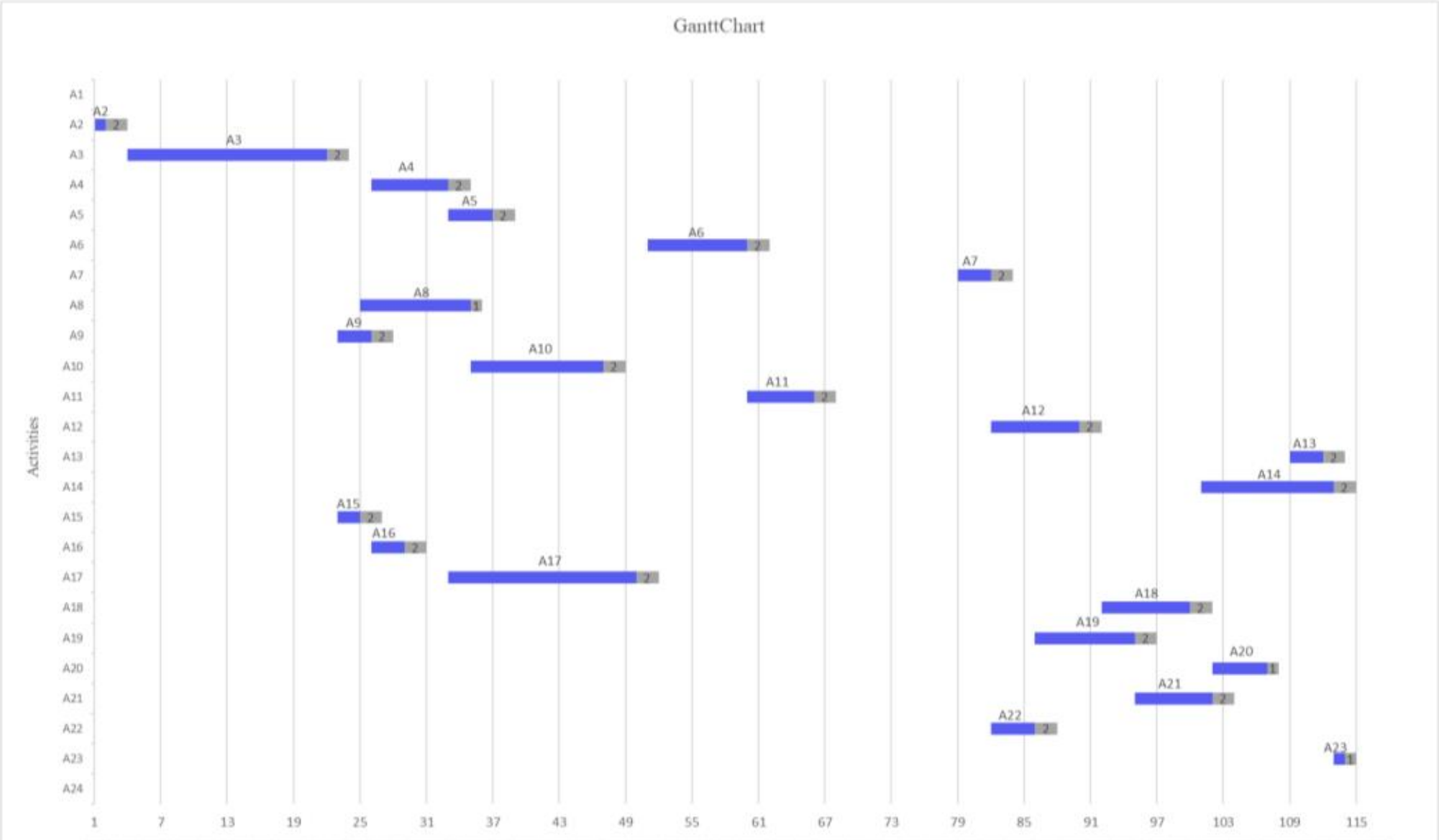

| | A24 | A23 | A22 | A21 | A20 | A19 | A18 | A17 | A16 | A15 | A14 | A13 | A12 | A11 | A10 | A9 | A8 | A7 | A6 | A5 | A4 | A3 | A2 | A1 |
|---|---|---|---|---|---|---|---|---|---|---|---|---|---|---|---|---|---|---|---|---|---|---|---|---|
| Start time | 114 | 113 | 82 | 95 | 102 | 86 | 92 | 33 | 26 | 23 | 101 | 109 | 82 | 60 | 35 | 23 | 25 | 79 | 51 | 33 | 26 | 4 | 1 | 1 |
| Duration | 0 | 1 | 4 | 7 | 5 | 9 | 8 | 17 | 3 | 2 | 12 | 3 | 8 | 6 | 12 | 3 | 10 | 3 | 9 | 4 | 7 | 18 | 1 | 0 |
| mode | 1 | 1 | 2 | 2 | 1 | 2 | 2 | 2 | 2 | 2 | 2 | 2 | 2 | 2 | 2 | 2 | 1 | 2 | 2 | 2 | 2 | 2 | 2 | 1 |
| finish time | 114 | 114 | 86 | 102 | 107 | 95 | 100 | 50 | 29 | 25 | 113 | 112 | 90 | 66 | 47 | 26 | 35 | 82 | 60 | 37 | 33 | 22 | 2 | 1 |

**Fig. 7.** The project Gantt chart

Similarly, the project Gantt chart is illustrated in Fig. 8A and Fig.8B. However, the difference between this Gantt chart and the last one in Fig. 7 is in *x*-axis representation. Meaning that, contrary to the last one, the project duration is shown in this figure in more detail. Moreover, four considered periods are depicted as black vertical lines on days 30, 60, 90, and 120. As such, amount of using different renewable and non-renewable resources, that is $BR\ (k=1)$, $BR\ (k=2)$, $WR\ (l=1)$, $WR\ (l=2)$, $WR\ (l=3)$, in addition to the cost of activities execution on each day are denoted in both Fig. 8A and Fig. 8B. It is worth noting that since the project is accomplished on day 114, resources are not required and consequently no cost is incurred. Therefore, there is no information from 114 to 120 in the table below the Gantt chart of Fig. 8B.

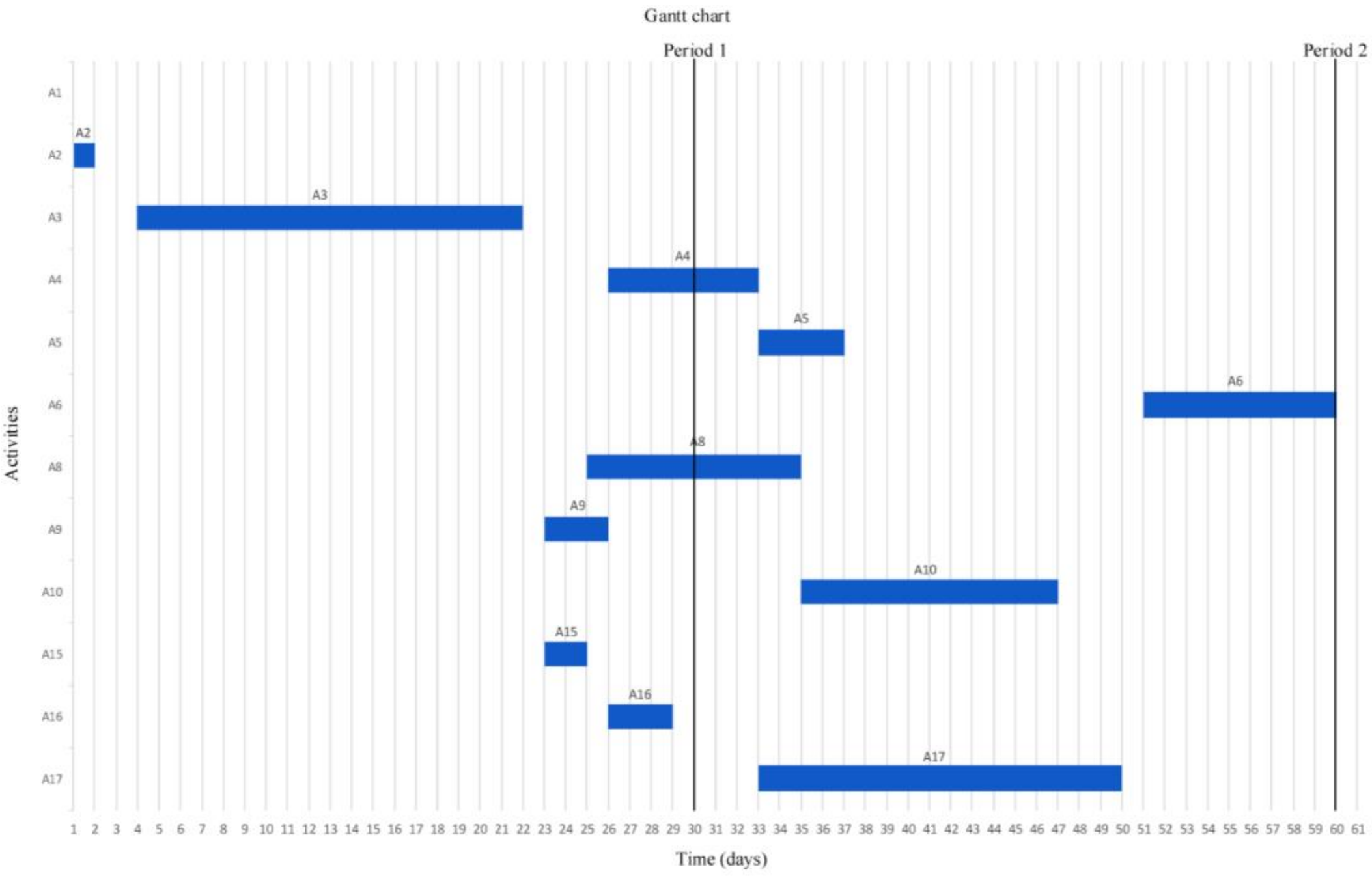

| | 1 | 2 | 3 | 4 | 5 | 6 | 7 | 8 | 9 | 10 | 11 | 12 | 13 | 14 | 15 | 16 | 17 | 18 | 19 | 20 | 21 | 22 | 23 | 24 | 25 | 26 | 27 | 28 | 29 | 30 |
|---|---|---|---|---|---|---|---|---|---|---|---|---|---|---|---|---|---|---|---|---|---|---|---|---|---|---|---|---|---|---|
| BR(k=1) | 3 | 3 | 3 | 6 | 6 | 6 | 6 | 6 | 6 | 6 | 6 | 6 | 6 | 6 | 6 | 6 | 6 | 6 | 6 | 6 | 6 | 6 | 3 | 3 | 11 | 14 | 14 | 14 | 14 | 14 |
| BR(k=2) | 2 | 2 | 2 | 5 | 5 | 5 | 5 | 5 | 5 | 5 | 5 | 5 | 5 | 5 | 5 | 5 | 5 | 5 | 5 | 5 | 5 | 5 | 4 | 4 | 13 | 17 | 17 | 17 | 17 | 17 |
| WR(l=1) | 6 | 6 | 6 | 8 | 8 | 8 | 8 | 8 | 8 | 8 | 8 | 8 | 8 | 8 | 8 | 8 | 8 | 8 | 8 | 8 | 8 | 8 | 17 | 17 | 30 | 25 | 25 | 25 | 25 | 25 |
| WR(l=2) | 8 | 8 | 8 | 10 | 10 | 10 | 10 | 10 | 10 | 10 | 10 | 10 | 10 | 10 | 10 | 10 | 10 | 10 | 10 | 10 | 10 | 10 | 16 | 16 | 31 | 29 | 29 | 29 | 29 | 29 |
| WR(l=3) | 7 | 7 | 7 | 9 | 9 | 9 | 9 | 9 | 9 | 9 | 9 | 9 | 9 | 9 | 9 | 9 | 9 | 9 | 9 | 9 | 9 | 9 | 15 | 15 | 29 | 27 | 27 | 27 | 27 | 27 |
| Cost | 3350 | 3350 | 3350 | 5450 | 5450 | 5450 | 5450 | 5450 | 5450 | 5450 | 5450 | 5450 | 5450 | 5450 | 5450 | 5450 | 5450 | 5450 | 5450 | 5450 | 5450 | 5450 | 6500 | 6500 | 1505 | 1610 | 1610 | 1610 | 1610 | 1610 |

| | 31 | 32 | 33 | 34 | 35 | 36 | 37 | 38 | 39 | 40 | 41 | 42 | 43 | 44 | 45 | 46 | 47 | 48 | 49 | 50 | 51 | 52 | 53 | 54 | 55 | 56 | 57 | 58 | 59 | 60 |
|---|---|---|---|---|---|---|---|---|---|---|---|---|---|---|---|---|---|---|---|---|---|---|---|---|---|---|---|---|---|---|
| BR(k=1) | 14 | 14 | 13 | 13 | 7 | 7 | 2 | 2 | 2 | 2 | 2 | 2 | 2 | 2 | 2 | 2 | 2 | 2 | 2 | 2 | 6 | 6 | 6 | 6 | 18 | 18 | 6 | 6 | 6 | 7 |
| BR(k=2) | 17 | 17 | 16 | 16 | 11 | 11 | 4 | 4 | 4 | 4 | 4 | 4 | 4 | 4 | 4 | 1 | 4 | 4 | 4 | 4 | 8 | 8 | 8 | 8 | 21 | 21 | 8 | 8 | 8 | 8 |
| WR(l=1) | 25 | 25 | 29 | 29 | 31 | 31 | 15 | 15 | 15 | 15 | 15 | 15 | 15 | 15 | 15 | 15 | 15 | 15 | 15 | 15 | 15 | 15 | 15 | 15 | 27 | 27 | 15 | 15 | 15 | 16 |
| WR(l=2) | 29 | 29 | 29 | 29 | 30 | 30 | 16 | 16 | 16 | 16 | 16 | 16 | 16 | 16 | 16 | 16 | 16 | 16 | 16 | 16 | 16 | 16 | 16 | 16 | 30 | 30 | 16 | 16 | 16 | 13 |
| WR(l=3) | 27 | 27 | 28 | 28 | 29 | 29 | 15 | 15 | 15 | 15 | 15 | 15 | 15 | 15 | 15 | 15 | 15 | 15 | 15 | 15 | 14 | 14 | 14 | 14 | 25 | 25 | 14 | 14 | 14 | 12 |
| Cost | 1610 | 1610 | 1595 | 1595 | 1360 | 1360 | 6200 | 6200 | 6200 | 6200 | 6200 | 6200 | 6200 | 6200 | 6200 | 6200 | 6200 | 6200 | 6200 | 6200 | 8050 | 8050 | 8050 | 8050 | 1800 | 1800 | 8050 | 8050 | 8050 | 7700 |

**Fig. 8A.** The daily resource usage and costs of activities during periods 1 and 2

Overall, the presented mathematical model takes into account multiple execution modes for performing activities. As a result, selecting the execution modes through which activities can be implemented faster, considering financial considerations affecting the project cash flow, leads to a reduction in project completion time. Moreover, calculating cash flow in each period enables more efficient resource allocation. Consequently, this approach ensures proper resource management throughout project execution, ultimately leading to greater project success and higher profits.

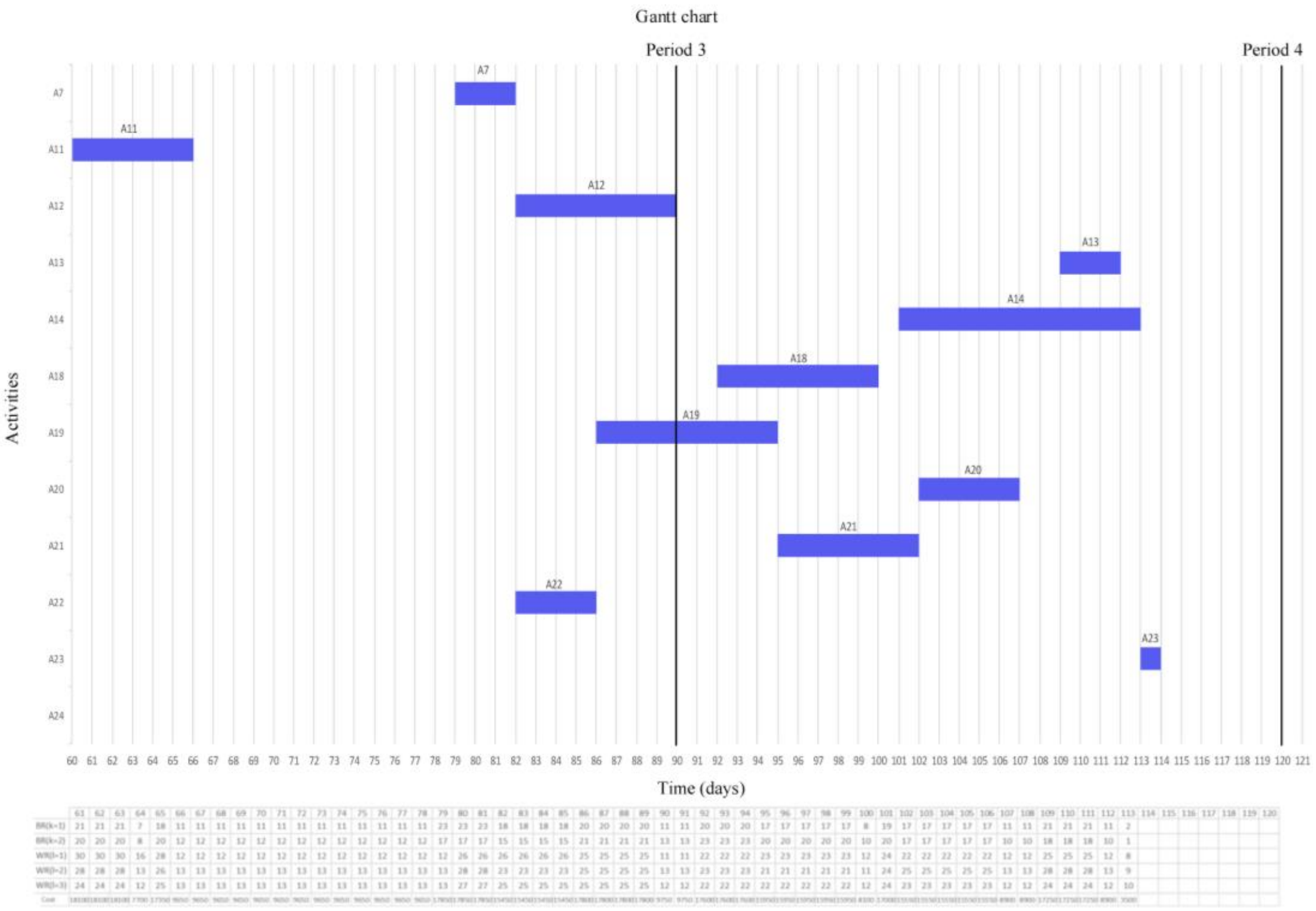

**Fig. 8B.** The daily resource usage and costs of activities during periods 3 and 4

### *5.3. Sensitivity Analysis*

To provide a more comprehensive insight, a sensitivity analysis is conducted herein. Equal weights of 0.5 are assigned to both objective functions. As such, the impact of different interest rates on the second objective function, final cash flow, is illustrated in Fig. 9. The considered interest rates on excess cash, delayed payments, long-term and short-term loans are 0.0125, 0.1, 0.06, and 0.075 respectively.

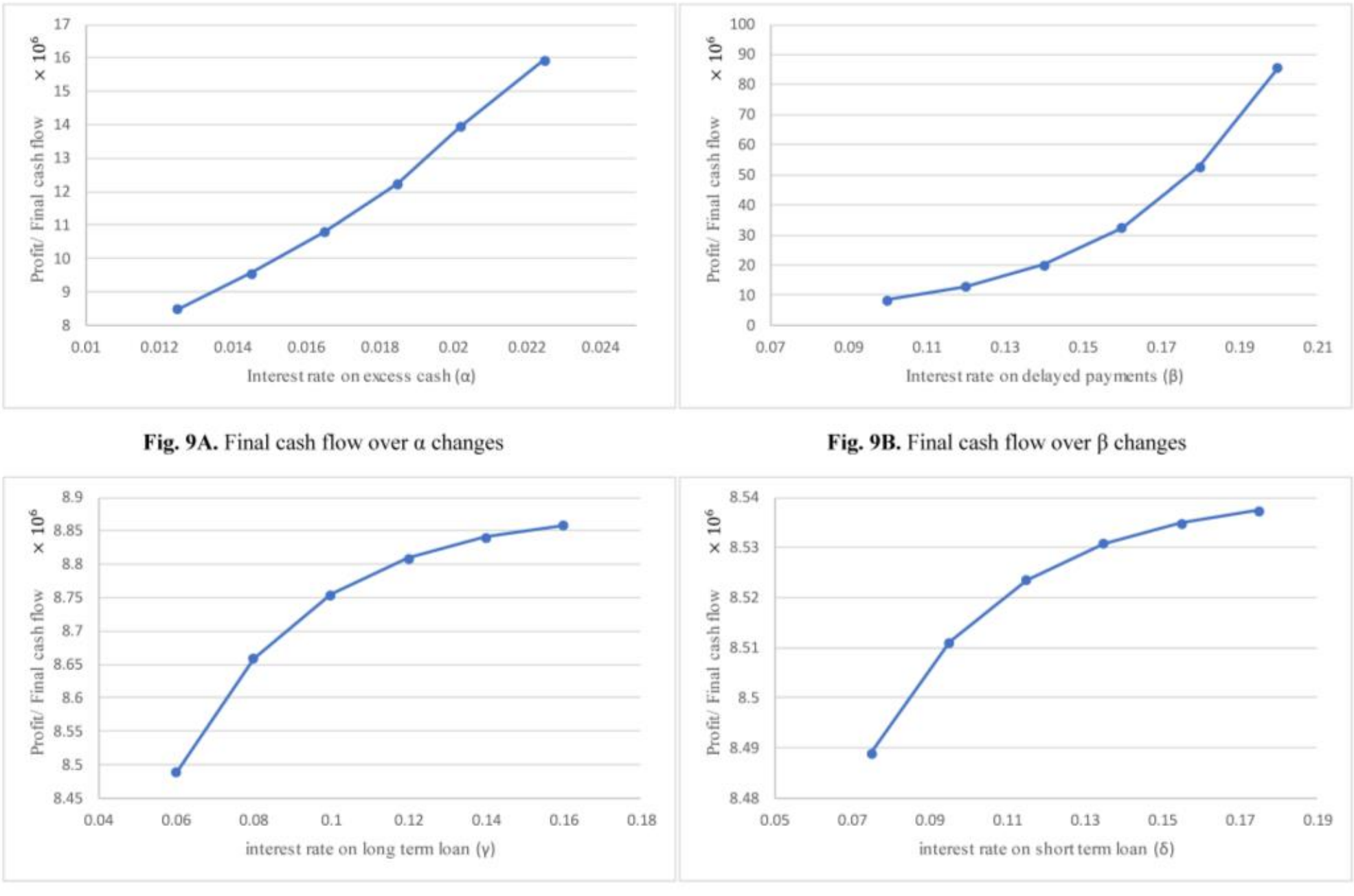

**Fig. 9A.** Final cash flow over α changes

**Fig. 9B.** Final cash flow over β changes

**Fig. 9C.** Final cash flow over γ changes

**Fig. 9D.** Final cash flow over δ changes

**Fig. 9.** Visual demonstration of sensitivity analysis of various interest rates

As it is shown in this figure, all final cash flows are increased by increasing interest rates. The final cash flow in Fig. 9A grows steadily by increasing interest rate on excess cash, which is started at 0.0125 and finished at 0.022. Likewise, it grows gradually from 0.1 to 0.16 in Fig. 9B. However, it escalates dramatically since then. Dissimilar to last figures, the final cash flows in Fig. 9C and Fig. 9D initially start increasing swiftly and subsequently grows slowly. Additionally, the corresponding information for various interest rates is presented in Table 5.

Table 5. Final cash flow under various interest rates

| $\alpha$ | Profit/ Final cash flow ($) | $\beta$ | Profit/ Final cash flow ($) | $\gamma$ | Profit/ Final cash flow ($) | $\delta$ | Profit/ Final cash flow ($) |
|---|---|---|---|---|---|---|---|
| 0.0125 | 8,488,946 | 0.1 | 848,8943 | 0.06 | 8,488,943 | 0.075 | 8,488,943 |
| 0.0145 | 9,557,932.85 | 0.12 | 12,854,620 | 0.08 | 8,659,288.74 | 0.095 | 851,1061.98 |
| 0.0165 | 10,798,640 | 0.14 | 20,186,530 | 0.1 | 8,755,180.2 | 0.115 | 8,523,611.83 |
| 0.0185 | 12,240,860 | 0.16 | 32,475,140 | 0.12 | 8,809,727.62 | 0.135 | 8,530,805.23 |
| 0.0202 | 13,939,040 | 0.18 | 52,728,880 | 0.14 | 8,841,071.38 | 0.155 | 8,534,969.1 |
| 0.0225 | 15,955,730 | 0.2 | 85,829,390 | 0.16 | 8,859,258.2 | 0.175 | 8,537,402.3 |

Performance of the peroposed mathematical model, considering the presented case study in this paper, is evaluated in comparison to the two recent models in cash flow problem presented by Alavipour and Arditi (2019) as well as Bruni and Hazir (2024) under α-level = 0.5. Table 6 compares the trade-off between time and profit in the mentioned studies and the current study. The main differences between this study and the research conducted by Alavipour and Arditi (2019) are considering initial capital, daily expenses, resource constraints, and multiple execution modes. Moreover, daily expenses, credit limit, periodic expenses, non-renewable resources, and multiple execution modes have not been addressed in Bruni and Hazir (2024). As shown in this table, in this paper, project is implemented approximately 10 and 15 days earlier than the studies conducted by Bruni and Hazir (2024) as well as Alavipour and Arditi (2019), respectively. Furthermore, profit of the project is nearly $ 5832805.125, which is nearly 41.71% and 49.41% more than the studies conducted by Bruni and Hazir (2024) as well as Alavipour and Arditi (2019), respectively.

Table 6. Comparative analysis of the proposed mathematical models in the current study, Bruni & Hazir (2024), and Alavipour (2019) under α-level = 0.5 in the implemented case study

| The current study | | Bruni and Hazir (2024) | | Alavipour (2019) | |
|---|---|---|---|---|---|
| Time (days) | Profit ($) | Time (days) | Profit ($) | Time (days) | Profit ($) |
| 67 | 5,217,831.025 | 83 | 4,052,095.265 | 93 | 3,844,303.207 |
| 72 | 5,281,017.476 | 88 | 4,083,006.360 | 94 | 3,856,606.436 |
| 78 | 5,335,091.622 | 91 | 4,098,461.907 | 95 | 3,878,909.664 |
| 83 | 5,355,188.460 | 96 | 4,115,965.572 | 99 | 3,887,663.509 |
| 88 | 5,363,227.134 | 99 | 4,153,543.371 | 102 | 3,937,848.795 |
| 94 | 8,444,475.034 | 101 | 4,192,845.980 | 107 | 4,018,443.162 |

Since multiple objective functions are considered in the proposed mixed-integer linear programming model, a novel approach, namely IVF-TH, is introduced to solve this multi-objective model. To demonstrate the efficiency of the proposed solution approach, a comparison is drawn between the new proposed approach and conventional fuzzy-TH approach. In this way, the two objective functions are investigated with these two approaches under different $\alpha$-levels, i.e. from 0.1 to 0.9, in Table7.

Table 7. Final cash flow ($) under different $\alpha$-levels across fuzzy- TH and IVF- TH

| α-level | Fuzzy- TH | IVF- TH | α-level | Fuzzy- TH | IVF- TH |
|---|---|---|---|---|---|
| 0.1 | 7,605,352.41 | 7,559,711.24 | 0.6 | 7,902,101.77 | 7,844,124.54 |
| 0.2 | 7,623,665.01 | 7,605,370.93 | 0.7 | 7,998,807.51 | 7,892,355.74 |
| 0.3 | 7,651,030.61 | 7,651,030.61 | 0.8 | 8,065,241.70 | 7,940,586.94 |
| 0.4 | 7,696,690.30 | 7,696,690.30 | 0.9 | 8,171,265.79 | 7,988,818.15 |
| 0.5 | 7,805,396.04 | 7,795,893.33 | | | |

More so, Fig. 10 confirm the better performance of the second objective function, which is maximizing the profit, for the IVF-TH approach in all α levels other than 0.3 and 0.4 in which both approaches have the same performance. Therefore, according to aforementioned tables and figures, the higher $\alpha$-levels are considered, the better performance of objective functions for both approaches can be observed. Nevertheless, compared to the Fuzzy-TH, the novel proposed approach in this paper, IVF-TH, has generally better performance for both objective functions in various α level.

To evaluate the performance of the proposed solution approaches, the most widely used RCPSP instances obtained from the PSPLIB (http://www.om-db.wi.tum.de/psplib/) and novel benchmark problem sets from MMLIB (http://www.projectmanagement.ugent.be) are utilized. The PSPLIB, from which a medium-size problem instance with 30 activities (i.e. J30) is considered, were created by ProGen (Kolisch and Sprecher 1997), and the MMLIB, from which large size benchmark problem sets with 50 and 100 activities (i.e. MM50 and MM100) are considered, were provided by Van Peteghem and Vanhoucke (2014). Table 8 compares the performance of IVF-TH to IVF-ABS, previously presented by Zarei et al. (2024), and a common multi-objective optimization method, namely the weighted sum method, considering IVF numbers in four $\alpha$-levels. The superiority of the proposed solution approach is clearly shown in both the objective function values as well as the CPU time. On average, the project of the mentioned case study was executed nearly four days and eight days earlier than IVF-ABS and IVF-Weighted sum, respectively. The second objective function (i.e. profit) is 3.50% and 8.42% more than IVF-ABS and IVF-Weighted sum, respectively. The CPU time of IVF-TH is 904 seconds, which is approximately 154 seconds and 244 seconds earlier than IVF-ABS and IVF-Weighted sum, respectively. Moreover, as the size of the problem instances increases, the performance of the IVF-TH approach is shown in this table.

Table 8. Comparative analysis of the IVF-TH, IVF-ABS, and IVF-Weighted sum under different $\alpha$-levels and problem size

| Problem instance | $\alpha$-levels | IVF-TH | | IVF-ABS | | IVF-Weighted sum | |
|---|---|---|---|---|---|---|---|
| | | Time (days) | Profit | Time (days) | Profit | Time (days) | Profit |
| Case study | $\alpha = 0.2$ | 67 | $ 7,605,370.93 | 71 | $ 7,484,917.703 | 75 | $ 6,955,955.110 |
| | $\alpha = 0.4$ | 67 | $ 7,696,690.30 | 71 | $ 7,461,219.413 | 75 | $ 7,132,552.406 |
| | $\alpha = 0.6$ | 67 | $ 7,844,124.54 | 70 | $ 7,500,935.568 | 74 | $ 7,229,888.181 |
| | $\alpha = 0.8$ | 66 | $ 7,940,586.94 | 70 | $ 7,588,309.966 | 74 | $ 7,323,443.462 |
| J30 | $\alpha = 0.2$ | 54 | 129,1399,116.25 | 57 | 1,080,289,137.66 | 60 | 921045246.88 |
| | $\alpha = 0.4$ | 52 | 1,310,682,577.47 | 57 | 1,090,856,341.45 | 58 | 1011829042.27 |
| | $\alpha = 0.6$ | 52 | 1,326,684,544.71 | 55 | 1112532327.63 | 58 | 1030686482.94 |
| | $\alpha = 0.8$ | 50 | 1,351,024,609.38 | 55 | 1131509411.22 | 58 | 1048730548.56 |
| MM50 | $\alpha = 0.2$ | 43 | 16,182,728,871.89 | 50 | 13,589,710,472.40 | 53 | 11,970,099,636.75 |
| | $\alpha = 0.4$ | 43 | 16,974,432,637.08 | 49 | 14854378480.92 | 53 | 12,097,932,254.14 |
| | $\alpha = 0.6$ | 41 | 19,422,174,331.39 | 49 | 14948251679.31 | 51 | 12,866,577,945.18 |
| | $\alpha = 0.8$ | 41 | 19,678,377,883.14 | 46 | 15273538032.58 | 47 | 13,364,370,937.76 |
| MM100 | $\alpha = 0.2$ | 61 | 44,792,107,291.52 | 65 | 41,308,462,136.50 | 73 | 40,111,269,354.17 |
| | $\alpha = 0.4$ | 60 | 45,816,099,081.04 | 65 | 43,604,430,038.00 | 72 | 41,196,670,077.72 |
| | $\alpha = 0.6$ | 59 | 46,689,151,785.00 | 64 | 44,160,610,876.17 | 72 | 42,496,430,055.70 |
| | $\alpha = 0.8$ | 58 | 47,182,356,655.70 | 63 | 44,915,995,626.50 | 70 | 43,111,720,396.45 |

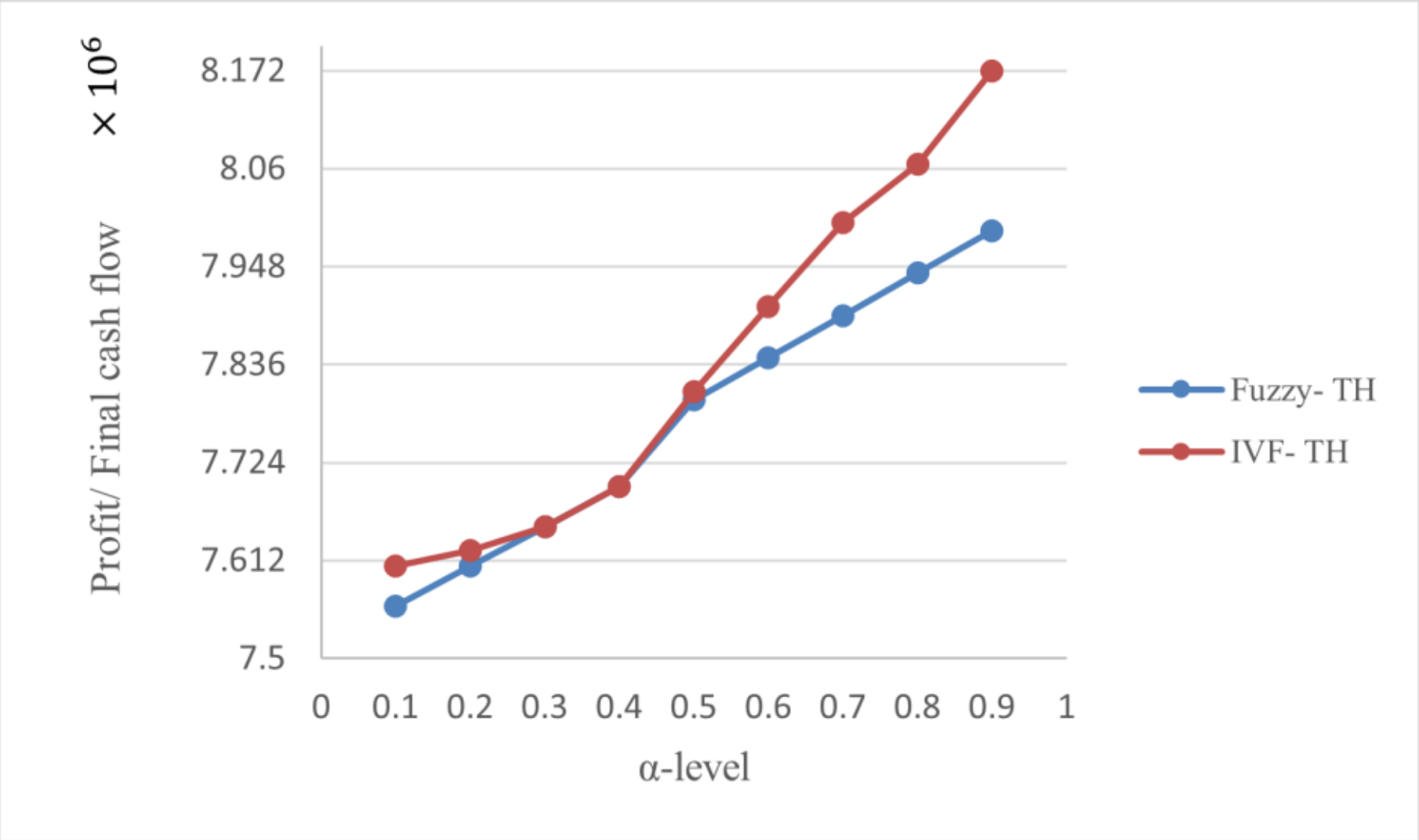

**Fig. 10.** Sensitivity analysis of final cash flow under various $\alpha$-levels

It is worth noting that, the presented bi-objective optimization model along with the proposed solution approach can schedule projects efficiently and calculate the cash flow of projects both daily and periodically. Given the relations between periods, this study can give rise to cost saving and have positive long-term impacts on profit of projects. In addition, utilizing the presented approach in this research, risks of resource shortages and consequently project halts are mitigated. In other words, proper project scheduling and project cash flow management in initial time periods and throughout the project life cycle can ensure the effective resource utilization in long-term.

## 6. Conclusions

As addressed in the literature, in spite of several attempts have been made, much more investigations need to be carried out in this regard. Therefore, a new multi-mode multi-objective mixed-integer linear programming model was proposed in this paper. As the proposed model was a multi-objective optimization model under fuzzy uncertainty, a new IVF-TH solution approach was introduced as well. A case study of construction projects in the oil and gas industry was presented to illustrate the proposed model applications. The model was solved in both fuzzy and crisp environments and the results were discussed. The results stemming from the model considering uncertainties were better than the crisp environment. Furthermore, the sensitivity analysis indicated that in comparison with the conventional fuzzy-TH approach, the proposed IVF-TH solution approach had generally better performance for both objective functions in various α levels. While the proposed model demonstrates significant potential in optimizing construction project scheduling, it is not without its limitations. Due to the prior experience of decision makers and managers, the probability distribution of some imprecise parameters might be available. In these circumstances, stochastic techniques can be employed to tackle the existing uncertainties. Moreover, inflation as a critical assumption that has considerable impacts on cash flow problems can be incorporated in the mathematical model to decrease the gap between theoretical and practical implications. While the presented case study demonstrated the applicability of the model in a specific area, extending its evaluation across diverse project types would provide valuable insights into its generalizability. Additionally, the current model primarily focuses on financial and resource constraints, leaving other real-world considerations unaddressed. Future research could address regulatory constraints to consider compliance with legal and safety standards, environmental factors to promote sustainability, and workforce limitations. Another practical assumption that can be taken into account in the extended version of the proposed model in future studies, is considering penalties for delayed repayment of loans affecting the results of the second objective function (profit).

## Data availability Statement

The authors confirm that the data supporting the findings of this study are available within the article.

**Acknowledgement**

The open-access publication fee for this research was supported by Deakin University.